\documentclass[a4paper,12pt]{article}  
\RequirePackage[OT1]{fontenc}  
\usepackage[natbib,noinfoline]{imsart}  
\usepackage[OT1]{fontenc}  
\usepackage[english]{babel}  
\usepackage[utf8]{inputenc}   
\usepackage{latexsym,amssymb,amsmath,epsfig,psfrag}  
  
\def\PP{\mathbb{P}} 
\def\RR{\mathbb{R}} 
\def\EE{\mathbb{E}} 
 
\def\ZZ{\mathbb{Z}}

\newcommand{\sfrac}[2]{\kern.1em
        \raise.5ex\hbox{$#1$}\kern-.1em
        /\kern-.15em\lower.25ex\hbox{$#2$}}

\def\bdes{\begin{description}}
\def\edes{\end{description}}

\newtheorem{defi}{Definition}[section]
\newtheorem{lemm}[defi]{Lemma}
\newtheorem{prop}[defi]{Proposition}

\newtheorem{theo}[defi]{Theorem}
\newtheorem{exem}{Example}
\newtheorem{rem}{Remark}
\newtheorem{conj}[defi]{Conjecture}
\newtheorem{question}[defi]{Question}
\newenvironment{dem}{\vskip 2mm\noindent {\it Proof} :}
                    {\hfill $\square$ \vskip 2mm \noindent} 
 
\newcommand{\eps}{\varepsilon}
\def\bex{\begin{exem} \em }
\def\eex{\end{exem} }
\def\brem{\begin{rem} \em}
\def\erem{\end{rem} }

\begin{document}

\begin{frontmatter}

\title{Submean variance bound for effective resistance of random electric networks} 

\runtitle{Variance bounds and random networks}
\thanks{Rapha\"el Rossignol was supported by the Swiss National Science
  Foundation grants 200021-1036251/1 and 200020-112316/1.} 
\author{\fnms{Itai} \snm{Benjamini}}
\address{The Weizmann Institute\\ Rehovot 76100\\  Isra\"el}
\author{\fnms{Rapha\"el} \snm{Rossignol} }
\address{Universit\'e de Neuch\^atel\\ Institut de Math\'ematiques\\ 11 rue Emile
  Argand\\ Case postale 158\\ 2009 Neuch\^atel, Suisse
}

\runauthor{I. Benjamini, R. Rossignol}

\begin{abstract}~: We study a model of random electric networks with Bernoulli
  resistances. In the case of the lattice $\ZZ^2$, we show that the point-to-point
  effective resistance between 0 and a vertex $v$ has a variance of order at most $(\log
  |v|)^{\frac{2}{3}}$ whereas its expected value is of order $\log |v|$, when
  $v$ goes to infinity. When $d\not =2$, expectation and variance are of the same
  order. Similar results are obtained in the context of $p$-resistance. The proofs rely on a modified Poincar\'e inequality due to
  Falik and Samorodnitsky \cite{FalikSamorodnitsky06}.
\end{abstract}


\end{frontmatter}

\section{Introduction}
The main goal of this paper is to study the effective resistance between two
finite sets of vertices in a random
electric network with i.i.d resistances. The infinite grid $\ZZ^d$ will be the
essential graph that we will focus on. Let us 
first briefly describe our notation for a deterministic electrical network (for more
background, see Doyle and Snell \cite{DoyleSnell84}, Peres
\cite{Peres97}, Lyons and Peres \cite{LyonsPeresinprogress} and Soardi \cite{Soardi94}). Let $G=(V,E)$ be
an unoriented  locally finite graph with an at most countable set of vertices $V$ and a set of
edges $E$ (we allow multiple edges between two vertices). Let
$r=(r_e)_{e\in E}$ be a collection of positive real numbers, which are called
resistances. To each edge $e$, one may associate two oriented edges, and we
shall denote by $\overrightarrow{E}$ the set of all these oriented edges. Let $A$ and $Z$ be
two finite, disjoint, non empty sets of vertices of $G$: $A$ will denote the source of the
network, and $Z$ the sink. A function $\theta$
on $\overrightarrow{E}$ is called a flow from $A$ to $Z$ with strength $\|\theta\|$ if it is
antisymmetric, i.e
$\theta_{\overrightarrow{xy}}=-\theta_{\overrightarrow{yx}}$, if it satisfies the
node law at each vertex $x$ of $V\setminus (A\cup Z)$:
$$\sum_{y\sim x}\theta_{\overrightarrow{xy}}=0\;,$$
and if the ``flow in'' at $A$ and the ``flow out'' at $Z$ equal $\|\theta\|$:
$$\|\theta\|=\sum_{a\in A}\sum_{\substack{y\sim a\\ y\not \in A}}\theta(\overrightarrow{ay})=\sum_{z\in Z}\sum_{\substack{y\sim z\\ y\not \in A}}\theta(\overrightarrow{yz})\;.$$
In this definition, it is assumed that all the vertices in $A$ are considered as a single one, as if they were linked by a wire with null resistance (and the same is true for $Z$). The effective resistance
$\mathcal{R}_r(A\leftrightarrow Z)$ may be defined in different ways, the
following is the most appropriate for us:
\begin{equation}
\label{eqdefresistance}
\mathcal{R}_r(A\leftrightarrow Z)=\inf_{\|\theta\|=1}\sum_{e\in E}
r_e\theta(e)^2\;,
\end{equation}
where the infimum is taken over all flows $\theta$ from $A$ to $Z$ with
strength 1. This infimum is always attained at what is called the
\emph{unit minimal (or wired)
  current} (see \cite{Soardi94} Theorem 3.25 p.~40). A current is a flow which
satisfies, in addition to the node law, Kirchhoff's loop law (see
\cite{Soardi94} p.~12). In finite graphs, currents are unique whereas in
infinite graphs, there may exist more than one current. But in $\ZZ^d$, for
instance, between two finite sets $A$ and $Z$, it is
known to be attained uniquely when the resistances are bounded
away from 0 and infinity (see Lyons and Peres \cite{LyonsPeresinprogress}
p.~82 or Soardi \cite{Soardi94} p.~39-43).

Electric networks have
been thoroughly studied by probabilists since there is a correspondence between
electrical networks on a given graph and reversible Markov chains on the same
graph. Let us introduce randomness on the electrical network itself by choosing the resistances independently and
identically distributed. That is to say, let $\nu$ be a probability measure on $\RR_+$,
and equip $\RR_+^E$ with the tensor product $\nu^{\otimes E}$. When the resistances are bounded away from 0 and $\infty$, it is easy to see that the mean of the
effective resistance is of the same order of that in the network where all
resistances are equal to 1. In fact, different realizations of this network
are ``roughly equivalent'' (see Lyons and Peres \cite{LyonsPeresinprogress}
p.~42), and for example, the associated random walks are of
the same type. Related results are those of Berger \cite{Berger02}, p.~550 and
Pemantle and Peres \cite{PemantlePeres96}, which give respectively sufficient conditions for almost sure recurrence of the network and a necessary and
sufficient condition for almost sure transience. In this paper, we are mainly 
concerned with the typical fluctuations of the function
$r\mapsto\mathcal{R}_r(A\leftrightarrow Z)$ around its mean when $A$ and $Z$
are ``far apart''. For simplicity, we choose to focus on the variance of the
effective resistance. Typically, we will take $A$ and $Z$ reduced to two
vertices far apart: $A=\{a\}$ and $Z=\{z\}$, and we shall note
$\mathcal{R}_r(a\leftrightarrow z)$ instead of
$\mathcal{R}_r(\{a\}\leftrightarrow \{z\})$. Following the terminology used in
First Passage Percolation (see \cite{Kesten84}), we shall call this the
\emph{point-to-point effective resistance} from $a$ to $z$.

In this paper, we prove that the type of fluctuations of the point-to-point effective resistance on
$\ZZ^d$ is
qualitatively different when $d=2$ and $d\not =2$. Indeed, when $d\not =2$, we
will see, 
quite easily, that these fluctuations are of the same order as its mean. On the other hand, when $G=\ZZ^2$, and the resistances are bounded away from 0 and $\infty$,
it is easy to show that the mean of $\mathcal{R}_r(0\leftrightarrow v)$ is of
order $\log |v|$, where $|v|$ stands for the $l^1$-norm of the vertex $v$ (see
section \ref{sec:Z2}). The main result of this paper is the following variance bound on $\ZZ^2$ when the resistances are
distributed according to a Bernoulli distribution bounded away from 0. 
\begin{theo}
\label{thmmain}
Suppose that $\nu=\frac{1}{2}\delta_a+\frac{1}{2}\delta_b$, with $0<a\leq
b<+\infty$. Let $E$ be the set of edges in $\ZZ^2$, and define
$\mu=\nu^{E}$. Then, as $v$ goes to infinity:
$$Var_\mu(\mathcal{R}_r(0\leftrightarrow v))=O\left((\log |v|)^{\frac{2}{3}}\right)\;.$$
\end{theo}

\noindent Here as in the rest of the article, when $f$ and $g$ are two functions on $\ZZ^d$, we use the notation ``$f(v)=O(g(v))$ as $v$ goes to infinity'' to mean there is a positive constant $C$ such that, for $v$ large enough,
$$f(v)\leq C g(v)\;.$$
We shall also use the notation ``$f(v)=\Theta(g(v))$'' to mean ``$f(v)=O(g(v))$ and $g(v)=O(f(v))$''.

The paper is organised as follows. In section \ref{secBernoulli} we introduce
the main tool of this paper, which is a modified Poincar\'e inequality due to
Falik and Samorodnitsky \cite{FalikSamorodnitsky06}. A first result is
given in Proposition \ref{propBKSelectric}, which announces our main result, on
$\ZZ^2$. Section \ref{sec:Z2} is devoted to the analysis of $\ZZ^2$: we prove our
main result, Theorem  \ref{thmmain} and
compare it to the simpler case of $\ZZ^d$, for $d\not = 2$. The choice of the
Bernoulli setting has been done for the sake of simplicity, but it is possible
to extend our variance bound to other distributions, and even to obtain the
corresponding exponential concentration inequalities. This is developped in
section \ref{sec:other}. In section \ref{sec:leftright}, we make some remarks
and conjectures on the a priori simpler case of the left-right resitance on the $n\times n$ grid. Finally,
section \ref{sec:penergy} is devoted to an extension of Theorem \ref{thmmain} to the non-linear
setting of $p$-networks.

\section{A general result in a Bernoulli setting}
\label{secBernoulli}
In this section, we suppose that $\nu$ is the Bernoulli probability measure
$\nu =\frac{1}{2}\delta_a+\frac{1}{2}\delta_b$, with $0<a\leq
b<+\infty$ and that for each collection of resistances $r$ in $\Omega=\{a,b\}^E$, there exists a unique  current
flow between two finite sets of vertices of the graph $G$. We want to bound from above the
variance of the effective resistance. Let us denote:
$$f(r)=\mathcal{R}_r(A\leftrightarrow Z)\;.$$
 A first idea is to use Poincar\'e
inequality, which in this setting is equivalent to Efron-Stein inequality
(see e.g. Steele \cite{Steele86} or An\'e et al. \cite{thesardsledoux}):
$$Var(f)\leq \sum_{e\in E}\left\|\Delta_e f\right\|_2^2\;,$$
where $\Delta_e$ is the following discrete gradient:
$$\Delta_e f(r)=\frac{1}{2}(f(r)-f(\sigma_er))\;,$$
$$(\sigma_er)_{e'}=\left\lbrace\begin{array}{ll}r_{e'}&\mbox{
    if } e'\not =e\\ b+a-r_e &\mbox{
    if }e'=e\end{array}\right.$$
Let $\theta_r$ be a flow attaining the minimum in the definition of
    $\mathcal{R}_r(a\leftrightarrow z)$. Using the definition of the effective
    resistance (\ref{eqdefresistance}),
\begin{eqnarray}
\nonumber f(\sigma_er)-f(r)&\leq &\sum_{e'\in
  E}(\sigma_er)_{e'}\theta_{r}(e')^2-\sum_{e'\in
    E}r_{e'}\theta_{r}(e')^2\;,\\
\label{eqmajodiscgrad}f(\sigma_er)-f(r)&\leq&(b-a)\theta_{r}(e)^2\;.
\end{eqnarray}
For any real number $h$, we denote by $h_+$ the number $\max\{h, 0\}$,
\begin{eqnarray}
\nonumber \sum_{e\in E}\left\|\Delta_e f\right\|_2^2&=&\frac{1}{2}\sum_{e\in E}\EE((f(\sigma_er)-f(r))_+^2)\;,\\
\label{eqmajoPoincare}\sum_{e\in E}\left\|\Delta_e f\right\|_2^2&\leq &\frac{(b-a)^2}{2}\EE(\sum_{e\in E}\theta_{r}(e)^4)\;.
\end{eqnarray}
It is quite possible that this last bound is sharp in numerous settings of
    interest, including $\ZZ^2$ (see section \ref{sec:Z2}), but in general we do
    not know how to evaluate the right-hand side of inequality
    (\ref{eqmajoPoincare}). We are just able to bound it from above using the fact that when $\theta$ is a unit current flow, $|\theta(e)|\leq 1$ for
every edge $e$. This last fact is intuitive, but for a formal proof, one can
    see Lyons and Peres \cite{LyonsPeresinprogress}, p.~49-50. Therefore
\begin{eqnarray}
\nonumber \sum_{e\in E}\left\|\Delta_e f\right\|_2^2&\leq &\frac{(b-a)^2}{2a}\EE(\sum_{e\in E}r_e\theta_{r}(e)^2)\;,\\
\label{eqmajoenergy} \sum_{e\in E}\left\|\Delta_e f\right\|_2^2&\leq
&\frac{(b-a)^2}{2a}\EE(f)\;.
\end{eqnarray}
 We have shown that the variance of
$f$ is at most of the order of its mean. It is possible to improve on this,
under some suitable assumption, by using the following
inequality, due to Falik and Samorodnitsky \cite{FalikSamorodnitsky06}:
\begin{equation}
\label{logsobFS}
Var(f)\log\frac{Var(f)}{\sum_{e\in E}\left\|\Delta_e f\right\|_1^2}\leq
2\sum_{e\in E}\left\|\Delta_e f\right\|_2^2\;.
\end{equation}
In order to state it as a bound on the variance of $f$, and to avoid
repetitions, we present it in a
slightly different way:
\begin{lemm}{\bf Falik and Samorodnitsky. }
\label{lemmFS}
Let $f$ belong to  $L^1(\{a,b\}^E)$. Suppose that
$\mathcal{E}_1(f)$ and $\mathcal{E}_2(f)$ are two real numbers such that:
$$\mathcal{E}_2(f)\geq \sum_{e\in E}\left\|\Delta_e f\right\|_2^2\;,$$
$$\mathcal{E}_1(f)\geq\sum_{e\in E}\left\|\Delta_e f\right\|_1^2\;,$$
and:
$$\frac{\mathcal{E}_2(f)}{\mathcal{E}_1(f)}\geq e\;.$$
Then,
$$Var(f)\leq
2\frac{\mathcal{E}_2(f)}{\log\frac{\mathcal{E}_2(f)}{\mathcal{E}_1(f)\log\frac{\mathcal{E}_2(f)}{\mathcal{E}_1(f)}}}\;.$$
\end{lemm}
\begin{dem}
Inequality (\ref{logsobFS}) is
proved by Falik and Samorodnitsky only for a finite set $E$, but it extends
straightforwardly to a countable set $E$, for functions in $L^1(\{a,b\}^E)$. Therefore, we have:
\begin{equation}
\label{eqmajobeforedisj}
Var(f)\log\frac{Var(f)}{\mathcal{E}_1(f)}\leq
2\mathcal{E}_2(f)\;.
\end{equation}
Now, consider the following disjunction:
\vspace*{2mm}

$\bullet$ either $Var(f)\leq \frac{\mathcal{E}_2(f)}{\log \frac{\mathcal{E}_2(f)}{\mathcal{E}_1(f)}}$,
\vspace*{2mm}

$\bullet$ or $Var(f)\geq \frac{\mathcal{E}_2(f)}{\log \frac{\mathcal{E}_2(f)}{\mathcal{E}_1(f)}}$, and plugging this inequality into (\ref{eqmajobeforedisj}) gives us:
$$Var(f)\leq 2\frac{\mathcal{E}_2(f)}{\log \frac{\mathcal{E}_2(f)}{\mathcal{E}_1(f)\log \frac{\mathcal{E}_2(f)}{\mathcal{E}_1(f)}}}\;.$$
\vspace*{2mm}

In any case, since $\mathcal{E}_2(f)/\mathcal{E}_1(f)\geq e$, the second
possibility is weaker than the first one, and we get:
$$Var(f)\leq 2\frac{\mathcal{E}_2(f)}{\log \frac{\mathcal{E}_2(f)}{\mathcal{E}_1(f)\log \frac{\mathcal{E}_2(f)}{\mathcal{E}_1(f)}}}\;.$$
\end{dem}
This inequality is very much in the spirit of an inequality by Talagrand
\cite{Talagrand94a} and could be called a modified Poincar\'e inequality (see
also 
\cite{BenaimRossignolarxiv06b}, \cite{BenaimRossignolarxiv06} and
\cite{Rossignol06} for more informations on such inequalities). The idea to use such a type of inequalities
in order to improve variance bounds is due to Benjamini, Kalai and Schramm
\cite{BenjaminiKalaiSchramm03} in
the context of First Passage Percolation. In our setting of random electric
networks, it allows us to show that as soon as the expected resistance is large but the
minimal energy flow via all but few resistances is small, then the variance of
the resistance is small compared to the expected resistance. This statement is
reflected in the following proposition, which is an introduction to the case
of $\ZZ^2$ in section \ref{sec:Z2}.
\begin{prop}
\label{propBKSelectric}
Let $G=(V,E)$ be
an unoriented  graph with an at most countable set of vertices $V$, a set of
edges $E$. Let $A$ and $Z$ be two disjoint non empty subsets of $V$. Let $a$ and $b$ be two positive real numbers, and $(r_e)_{e\in E}$
be i.i.d resistances with common law $\frac{1}{2}\delta_a+\frac{1}{2}\delta_b$.  Let $E_m$ be any subset of $E$ such that $E_m^c$ is finite. Define:
$$\alpha_m=\sup_{e\in E_m}\EE(r_e\theta_r^2)\;,$$
$$\beta_m=\frac{|E_m^c|}{\EE(\mathcal{R}_r(A\leftrightarrow Z))}\;,$$
and
$$\eps_m=\left(\frac{b-a}{a}\right)^2\alpha_m+(b-a)^2\beta_m\;.$$
Suppose that $\eps_m<1$. Then,
$$Var(\mathcal{R}_r(A\leftrightarrow Z))\leq 2K\frac{\EE(f)}{\log \frac{K}{\eps_m\log \frac{K}{\eps_m}}}\;,$$
where $K=\sup\left\lbrace\frac{(b-a)^2}{2a},e\right\rbrace$.
\end{prop}
\begin{dem}
We want to use Lemma \ref{lemmFS}. Define:
$$f(r)=\mathcal{R}_r(a\leftrightarrow z)\;.$$
Let us evaluate the terms $\sum_{e\in E}\left\|\Delta_e
  f\right\|_1^2$ and $\sum_{e\in E}\left\|\Delta_e
  f\right\|_2^2$. We have already seen in inequality (\ref{eqmajoenergy}) that:
$$\sum_{e\in E}\left\|\Delta_e
  f\right\|_2^2\leq \frac{(b-a)^2}{2a}\EE(f)\;,$$
and so, define:
$$\mathcal{E}_2(f)=K\EE(f)\;,$$
where $K=\sup\left\lbrace\frac{(b-a)^2}{2a},e\right\rbrace$. Besides,
\begin{eqnarray*}
\sum_{e\in E} \left\|\Delta_e f\right\|_1^2&=&\sum_{e\in E_m} \left\|\Delta_e
  f\right\|_1^2+\sum_{e\in E_m^c} \left\|\Delta_e f\right\|_1^2\;,\\
&\leq&\frac{b-a}{a}\sum_{e\in E}\left\|\Delta_e f\right\|_1\sup_{e\in
  E_m}\EE(r_e\theta_r^2)+|E_m^c|\sup_{e\in E}\left\|\Delta_e
  f\right\|_1^2\;,\\
&=&\frac{b-a}{a}\alpha_m\sum_{e\in E}\left\|\Delta_e f\right\|_1+\sup_{e\in E}\left\|\Delta_e
  f\right\|_1^2\beta_m\EE(f)\;,
\end{eqnarray*}
Recall that when $\theta$ is a unit current flow, $|\theta(e)|\leq 1$ for
every edge $e$. Therefore,
$$\sup_{e\in E}\left\|\Delta_e f\right\|_1 \leq (b-a)\;.$$
Also, using inequality (\ref{eqmajodiscgrad}),
\begin{eqnarray*}
\sum_{e\in E}\left\|\Delta_e f\right\|_1&=&\sum_{e\in E}\EE((f(\sigma_er)-f(r))_+)\;,\\
&\leq &(b-a)\EE(\sum_{e\in E}\theta_{r}(e)^2)\;,\\
&\leq &\frac{b-a}{a}\EE(f)\;.\\
\end{eqnarray*}
Therefore,
$$\sum_{e\in E} \left\|\Delta_e f\right\|_1^2\leq \EE(f)\eps_m\;.$$
Define 
$$\mathcal{E}_1(f)=\EE(f)\eps_m\;.$$
The assumption $\eps_m<1$ ensures that $\mathcal{E}_2(f)/\mathcal{E}_1(f)\geq
e$. We conclude by applying Lemma \ref{lemmFS}.
\end{dem}
Unfortunately, it is not very easy to bound the terms $\alpha_m$ and $\beta_m$
in an efficient way, essentially because in a random setting, we have no good
bound on the amount of current through a particular edge. For example, in section \ref{sec:Z2}, we will have to
resort to an averaging trick, and we shall not be able to use directly
Proposition \ref{propBKSelectric}. Nevertheless, for some interesting graphs such as trees, but also the lattices $\ZZ^d$, exact calculations are available when all
resistances are equal. Therefore, Proposition \ref{propBKSelectric} 
would become more helpful if one could prove the following stability result,
which we deliberately state in an informal way.
\begin{question}
Assume the flow on the fixed
resistance 1 environment on a graph satisfies the condition that, except for a small set of edges, only $o(1)$ flow goes via all the
other edges, then is the same true for a perturbed environment ? 
\end{question}

\section{The $\ZZ^d$ case}
\label{sec:Z2}

It is natural to inspect the Bernoulli setting on the most studied electrical
networks, which are $\ZZ^d$, $d\geq 1$. Here, we focus on the point-to-point
resistance between the origin and a vertex $v$ when $v$ goes to infinity
($+\infty$ when $d=1$). Let us
denote it as $f_v$:
$$f_v(r)=\mathcal{R}_r(0\leftrightarrow v)\;.$$

\subsection{The $\ZZ^d$ case for $d\not = 2$}

When $d\not = 2$, one can see
easily that the variance of $f_v$ is of the same order as its mean
value (when $b>a$). Indeed, when $d=1$, $f_n$ is just $na+(b-a)B_n$, where $B_n$ is a
random variable of binomial distribution with parameters $n$ and $1/2$.

When
$d\geq 3$, remark first that, denoting $\overline{a}=(a,a,\ldots)$,
\begin{equation}
\label{eqcompresist}
a\mathcal{R}_{\overline{1}}(0\leftrightarrow v)\leq\mathcal{R}_r(0\leftrightarrow v)\leq
b\mathcal{R}_{\overline{1}}(0\leftrightarrow v)\;.
\end{equation}
Therefore, the mean of $f_v$ is of the same order (up to a multiplicative
constant) as in the network where all resistances equal 1. Thus, when
$d\geq 3$, the mean of $f_v$ is of order $\Theta(1)$ (see Lyons and Peres
\cite{LyonsPeresinprogress} p.~39-40). The variance of $f_v$ is also of order
$\Theta(1)$ when $b>a$, as follows from the following simple lemma.
\begin{lemm}
\label{lemminovariance}
Let $G=(V,E)$ be a unoriented connected graph, $s\in V$ a vertex with finite degree $D$ and $Z$ be a finite subset of
$V$ such that $s\not\in Z$. If the resistances on $E$ are independently
distributed according to a symmetric Bernoulli law on $\{a,b\}$, with $b\geq a>0$,
$$Var(\mathcal{R}_r(s\leftrightarrow Z))\geq C(b-a)^2\;,$$
where $C$ is a positive constant depending only on $D$.
\end{lemm}
\begin{dem}
Denote by $\mathcal{D}={e_1,\ldots ,e_D}$ the $D$ edges incident to $s$. For any $r$ in
$\RR_+^E$, denote 
$$f(r)=\mathcal{R}_r(s\leftrightarrow Z)\;,$$
and let $(b^{(D)},r^{-D})$ be the the set of resistances obtained from $r$ by
switching all resistances on $\mathcal{D}$ to $b$. One has:
\begin{eqnarray*}
Var(f)&\geq &\EE((f-\int f\;dr_{e_1}\ldots dr_{e_D})^2)\;,\\
&\geq &\frac{1}{2^{D}}\EE((f(b^{(D)},r^{-D})-\int
f\;dr_{e_1}\ldots dr_{e_D})^2)\;,\\
&\geq&\left(\frac{1}{2^{D}}\right)^2\EE((f(b^{(D)},r^{-D})-(f(a^{(D)},r^{-D}))^2)\;,\\
&\geq&\left(\frac{1}{2^{D}}\right)^2(b-a)^2\EE((\sum_{e\in\mathcal{D}}\theta_{(b^{(D)},r^{-D})}(e)^2)^2)\;,\\
&\geq&\left(\frac{1}{2^{D}}\right)^2(b-a)^2\frac{1}{D^2}\;.
\end{eqnarray*}
\end{dem}

\subsection{The case of $\ZZ^2$: some heuristics}
Now, let us examine the case of $\ZZ^2$. When the resistances are bounded away from 0 and infinity, $f_v$, and therefore its
expectation, is of order $\Theta(\log |v|)$. Indeed, equation
(\ref{eqcompresist}) implies that it is of the same order as the resistance on
the network where all resistances equal 1. This more simple resistance can be
explicitly computed using Fourier transform on the lattice $\ZZ^2$ (see Soardi
\cite{Soardi94} p.~104-107). In a more simple way, it can be easily bounded from
below by using Nash-Williams inequality, and
from above by embedding a suitable tree in $\ZZ^2$ (see Doyle and
Snell \cite{DoyleSnell84} p. 85, or alternatively Lyons and Peres
\cite{LyonsPeresinprogress} p.~39-40). A more complicated question to address is the
existence of a precise limit of the ratio $\EE(f_v)/\log |v|$. This would lead
to an analog of the ``time constant'' arising in the context of First
Passage Percolation (see Kesten \cite{Kesten84}). Closely related questions
are the existence of an asymptotic shape and, if it exists, whether it is an euclidean ball or not. We believe that the time
constant and the asymptotic shape exist.
\begin{conj}
Define 
$$B_t=\{v\in\ZZ^2\mbox{ s.t. }\mathcal{R}_r(0\leftrightarrow v)\leq
t\}\;.$$
There exists a non empty, compact subset of $\RR^2$, $B_0$ such that, for every positive
number $\eps$,
$$(1-\eps)B_0\subset \frac{1}{\log t}B_t\subset (1+\eps)B_0\;.$$
\end{conj}

What about the order of the variance of
$f_v$, when $v$ goes to infinity ? Reasonably, it should be of order $\Theta(1)$. Since we did not manage to prove this, we state it as a conjecture:
\begin{conj}
Suppose that $\nu=\frac{1}{2}\delta_a+\frac{1}{2}\delta_b$, with $0<a\leq
b<+\infty$. Let $E$ be the set of edges in $\ZZ^2$, and define
$\mu=\nu^{E}$. Then, as $v$ tends to infinity:
$$Var_\mu(\mathcal{R}_r(0\leftrightarrow v))=\Theta(1)\;.$$
\end{conj}

A first intuitive support to this conjecture comes from inequality
(\ref{eqmajoPoincare}).  It is quite possible that it gives a bound of order
$O(1)$.  Indeed, this would be the case if the current in the perturbed
environment remained ``close'' (for example at a $l^4$-distance of order 1) to
the current in the uniform network (with all resistances equal to 1).

A second support to this conjecture comes from the analysis of the graph
$G_n=(V_n, E_n)$. This one arises when one applies in a classical way the Nash-Williams inequality to
get a lower bound on the resistance between the origin and the border of the
box $\{-n,\ldots,n\}\times\{-n,\ldots,n\}$. The
set of vertices $V_n$ is just $\{0,\ldots, n\}$. For $i$ in $\{0,\ldots,
n-1\}$, draw $2i+1$ parallel edges between $i$ and $i+1$, and call them
$e_{i,1},\ldots, e_{i,2i+1}$. This is a Parallel-Series electric
network, and the effective resistance is easy to compute:
$$\mathcal{R}_r(0\leftrightarrow
n)=\sum_{i=0}^{n-1}\frac{1}{\sum_{k=1}^{2i+1}\frac{1}{r_{i,k}}}\;.$$
One can show the following result.
\begin{prop}
\label{propparser}
If the resistances on $G_n$ are independently
distributed according to a symmetric Bernoulli law on $\{a,b\}$, with $b\geq
a>0$,
$$\EE(\mathcal{R}_r(0\leftrightarrow n))=\Theta(\log n)\;,$$
and
$$Var(\mathcal{R}_r(0\leftrightarrow n))=\Theta(1)\;.$$
\end{prop}
\begin{dem}
To shorten the notations, we treat the case $a=1/2$ and $b=1$. For any $r$ in
$\{a,b\}^E$, denote 
$$f(r)=\mathcal{R}_r(0\leftrightarrow n)\;.$$
The estimate on the mean is obvious. The estimate on the variance is easy
too. First note that
$$Var(f)=\sum_{i=0}^{n-1}Var\left(\frac{1}{\sum_{k=1}^{2i+1}\frac{1}{r_{i,k}}}\right)\;.$$
Denote by $Y_i$ the random variable:
$$Y_i=\frac{1}{\sum_{k=1}^{2i+1}\frac{1}{r_{i,k}}}\;.$$
Remark that:
$$\sum_{k=1}^{2i+1}\frac{1}{r_{i,k}}=2i+1+B_{2i+1}\;,$$
where $B_{2i+1}$ has a binomial distribution with parameters $2i+1$ and
$1/2$. Therefore, denoting 
$$N_i=\frac{B_{2i+1}-\frac{2i+1}{2}}{\sqrt{\frac{2i+1}{4}}}\;,$$
when $i$ tends to infinity, $N_i$ converges weakly to a standard Gaussian
variable, and:
$$Y_i=\frac{1}{3\frac{(2i+1)}{2}+\sqrt{\frac{2i+1}{4}}N_i}\;.$$
Therefore,
$$Y_i=\frac{1}{9(i+1/2)^{3/2}}\left\lbrack3\sqrt{2i+1}\left(\frac{1}{1+\frac{N_i}{3\sqrt{2i+1}}}-1\right)\right\rbrack+\frac{1}{3(i+1/2)}\;.$$
Define:
$$Z_i=3\sqrt{2i+1}\left(\frac{1}{1+\frac{N_i}{3\sqrt{2i+1}}}-1\right)\;.$$
Since the sequence $(N_i)$ is weakly convergent, it is bounded in
probability. Hence, using that:
$$\frac{1}{1+x}-1+x=O(x^2)\;,$$
as $x$ goes to zero, we deduce that $\sqrt{2i+1}(Z_i+N_i)$ is bounded in
probability, and therefore $Z_i$ converges in distribution to a standard
Gaussian random variable. Using the concentration properties of the binomial distribution, it is easy to
show that $Z_i$ and $Z_i^2$ are asymptotically uniformly integrable. This
implies that the variance of $Z_i$ tends to 1 as $i$ tends to infinity. Thus,
$$Var(Y_i)=\frac{1}{9^2(i+1/2)^{3}}Var(Z_i)=\frac{1}{9(i+1/2)^{3}}(1+o(1))\;,$$
and consequently,
$$Var(f)=\Theta(1)\;.$$
\end{dem}
Of course, on $\ZZ^2$, things are more difficult to compute.

\subsection{The case of $\ZZ^2$: proof of Theorem \ref{thmmain}}
We shall prove below that the variance of $f_v$ is of order $O((\log
|v|)^{\frac{2}{3}})$. We shall proceed very much as in
\cite{BenjaminiKalaiSchramm03}, resorting to an averaging trick to trade the study of
$f_v$ against the study of a randomized version of it.

\vskip 2mm\noindent {\it Proof of Theorem \ref{thmmain}} :
Let $m$ be a positive integer to be fixed later, and $z$ a random variable, independent from the edge-resistances, and
distributed according to $\mu_m$, the uniform distribution on the box $B_m=[0,m-1]^2\cap \ZZ^2$. Define 
$$
\tilde f(z,r):= \mathcal{R}_r(z,v+z)\,.
$$
We think of $\tilde f$ as a function on the
space $\tilde\Omega:=B_m\times\{a,b\}^{E}$ which is endowed with the
probability measure $\mu_m\otimes \nu^{\otimes E}$. The first thing to note is that $f_v$
and $\tilde f$ are not too far appart. To see this, we can use the following triangle
inequality (see \cite{LyonsPeresinprogress}, exercise~2.65 p.~67), which holds for every
three vertices $x,y,z$ in $\ZZ^2$, and any $r\in\Omega$,
\begin{equation}
\label{eqpseudotriangular}
\mathcal{R}_r(x,z)\leq
\mathcal{R}_r(x,y)+\mathcal{R}_r(y,z)\;.
\end{equation}
Therefore, taking $L^2$-norms in $L^2(\mu_m\otimes\nu^{\otimes E})$, and
noting that $|z|\leq 2m$,
\begin{eqnarray*}
\|f_v-\tilde f\|_2&\leq& \|\mathcal{R}_r(0,z)\|_2+\|\mathcal{R}_r(z,z+v)\|_2\;,\\
&\leq&C\log m\;,
\end{eqnarray*}
where $C$ is a universal constant. Noting that $\tilde f$ has the same
expectation as $f_v$, we get thus:
\begin{eqnarray*}
\|f_v-\EE(f_v)\|_2&\leq &\|f_v-\tilde f\|_2+\|\tilde f-\EE(\tilde f)\|_2\;,\\
&\leq& C\log m+\|\tilde f-\EE(\tilde f)\|_2\;.
\end{eqnarray*}
Therefore:
\begin{equation}
\label{eqmajofftilde}
Var(f_v)\leq (C\log m)^2+2C\log m\sqrt{Var(\tilde f)}+Var(\tilde f)\;.
\end{equation}

Now, we want to bound the variance of $\tilde f$ from above.  Define:
$$\EE_{\mu}(\tilde f)=\int f(z,r)\;d\mu_m(z)\;,$$
$$\EE_{\nu}(\tilde f)=\int f(z,r)\;d\nu^{\otimes E}(r)\;,$$
$$Var_{\mu}(\tilde f)=\EE_\mu (f(z,r)^2)-\EE_\mu(f(z,r))^2\;,$$
$$Var_{\mu}(\tilde f)=\EE_\nu (f(z,r)^2)-\EE_\nu(f(z,r))^2\;.$$
Then, we split the
variance of $\tilde f$ into two parts: the one due to $z$ and the other due to $r$.
\begin{equation}
\label{eqdecompvar}
Var(\tilde f)=\EE_\nu(Var_\mu(\tilde f))+Var_{\nu}(\EE_{\mu}(\tilde f))\;.
\end{equation}
Thanks to the triangle inequality (\ref{eqpseudotriangular}),
\begin{equation}
\label{eqmajovarz}
\EE_\nu(Var_\mu(\tilde f))\leq (C\log m)^2\;.
\end{equation}
To bound the last term of the sum in (\ref{eqdecompvar}), we apply Lemma
\ref{lemmFS} to $\EE_{\mu}(\tilde f)$. Remark that, thanks to Jensen's
inequality,
$$\left\|\Delta_e \EE_{\mu}(\tilde f)\right\|_1\leq \left\|\Delta_e\tilde f
\right\|_1\;,$$
where the first $L^1$-norm integrates against $\nu^{\otimes E}$, and the
second one integrates against $\mu_m\otimes\nu^{\otimes E}$. Also,
$$\left\|\Delta_e \EE_{\mu}(\tilde f)\right\|_2\leq \left\|\Delta_e\tilde f
\right\|_2\;.$$
Let us denote by 
$\theta_r^z$ the unit current flow from $z$ to $z+v$, when the resistances are
$r$. Using inequality (\ref{eqmajodiscgrad}), and the translation invariance
of this setting, we get, for every edge $e$:
\begin{eqnarray*}
\left\|\Delta_e \tilde f\right\|_1 & \leq & (b-a)\EE(\theta_r^z(e)^2)\;,\\
& = & (b-a)\EE(\theta_r^0(e-z)^2)\;,\\
& =&\frac{b-a}{(m+1)^2}\sum_{z_0\in B_m}\EE(\theta_r(e-z_0)^2)\;.
\end{eqnarray*}
Now we claim that:
\begin{equation}
\label{eqmajoenergybox}
\forall e\in E,\;\sum_{z_0\in B_m}\EE(\theta_r(e-z_0)^2)\leq \frac{5b(m+1)}{a}\;.
\end{equation}
Assuming this claim, we have:
\begin{equation}
\label{eqmajonorm1}
\sup_{e\in E}\left\|\Delta_e \tilde f\right\|_1=O\left(\frac{1}{m}\right)\;.
\end{equation}
To see that claim (\ref{eqmajoenergybox}) is true, let $e_-$ be the left-most
or lower end-point of $e$ and remark that the set of
edges described by $e-z$ is included in the set of edges of the box
$B_m^e:=e_-+B_{m+1}$. Therefore,
$$\sum_{z_0\in B_m}\EE(\theta_r(e-z_0)^2)\leq \sum_{e'\subset B_{m}^e}\EE(\theta_r(e')^2)\;.$$
Let $\partial B_m^e$ be the (inner) border of the box $B_m^e$:
$$\partial B_m=\{0\}\times [0,m]\cup\{m\}\times [0,m]\cup [0,m]\times\{0\}
\cup [0,m]\times\{m\}\;.,$$
$$\partial B_m^e=e_-+\partial B_m\;.$$
First, suppose that neither $0$ nor $v$ belongs to $B_m^e$. We define a flow
$\eta$ from $0$ to $v$ such that:
$$\eta(e')=\theta_r(e')\mbox{ if }e'\not\subset B_m^e\;,$$
$$\eta(e')=0\mbox{ if }e'\subset B_m^e\mbox{ and }e'\not\subset \partial B_m^e\;.$$
These conditions do not suffice to determine uniquely the flow $\eta$, but one
can then choose the flow on $\partial B_m^e$ that minimizes the energy
$\sum_{e'\in\partial B_m^e}r_{e'}\eta(e')^2$. This is the current flow on
$\partial B_m^e$ when the flow entering and going outside 
$\partial B_m^e$ is fixed by $\theta_r$. For a formal proof of the existence
of such a flow, see Soardi \cite{Soardi94} Theorem 2.2, p.~22. This flow has a strength less than 1. Therefore, the flow through each edge of
$\partial B_m^e$ is less than 1, and:
\begin{eqnarray*}
\sum_{e'\subset B_{m}^e}r_e\theta_r(e)^2&=& \sum_{e\in
  E}r_e\theta_r(e)^2-\sum_{e\not\subset B_m^e}r_e\theta_r(e)^2\;,\\
&\leq& \sum_{e\in
  E}r_e\eta(e)^2-\sum_{e\not\subset B_m^e}r_e\theta_r(e)^2\;,\\
&=&\sum_{e\subset\partial B_m^e}r_e\eta(e)^2\;,\\
&\leq&4b(m+1)\;.
\end{eqnarray*}
Therefore,
$$\sum_{e'\subset B_{m}^e}\theta_r(e)^2\leq \frac{4b(m+1)}{a}\;.$$

Now, suppose that $a$ or $v$ belongs to $B_m^e$. Let us say $v$ belongs to
$B_m^e$, the other situation being symmetrical. We define a flow
$\eta$ from $0$ to $v$ as before,
except that instead of assigning 0 to each value inside $B_m^e$, we keep a
path from $\partial B_m^e$ to $v$ on which the flow is assigned to 1, directed
towards $v$. The same considerations as before lead to.
$$\sum_{e'\subset B_{m}^e}\theta_r(e)^2\leq \frac{5b(m+1)}{a}\;.$$
And claim (\ref{eqmajoenergybox}) is proved.

Finally notice that, using inequality (\ref{eqmajodiscgrad}):
$$\sum_{e\in E}\left\|\Delta_e \tilde f\right\|_1\leq\frac{b-a}{a}\EE(\tilde
f)\;,$$
$$\sum_{e\in E}\left\|\Delta_e \tilde f\right\|_2^2\leq\frac{(b-a)^2}{2a}\EE(\tilde
f)\;.$$
Recall that 
$$\EE(\tilde f)=\EE(f_v)=\Theta(\log |v|)\;.$$
Therefore, there exist constants $K$  and $K'$ such that:
\begin{eqnarray*}
\sum_{e\in E}\left\|\Delta_e
  \tilde f\right\|_1^2&\leq &\sum_{e\in E}\left\|\Delta_e
  \tilde f\right\|_1.\sup_{e\in E}\left\|\Delta_e \tilde f\right\|_1\;,\\
&\leq &K\frac{\log |v|}{m}\;,
\end{eqnarray*}
and:
$$\sum_{e\in
  E}\left\|\Delta_e \tilde f\right\|_2^2\leq K'\log |v|\;.$$
Denoting $\mathcal{E}_1(\tilde f)=K\frac{\log |v|}{m}$ and
  $\mathcal{E}_2(\tilde f)=K'\log |v|$, the hypotheses of Lemma \ref{lemmFS} are
  fulfilled, at least for $m$ larger than $e\frac{K}{K'}$. Assume that $m$ is
  a function of $|v|$ which goes to infinity when $|v|$ goes to infinity. Lemma \ref{lemmFS} together with inequalities (\ref{eqdecompvar})
  and (\ref{eqmajovarz}) gives us:
$$Var (\tilde f)\leq  O\left(\frac{\log |v|}{\log m}+(\log m)^2\right)\;.$$
Thus, choose $m$ the greatest integer such that $\log m \leq (\log |v|)^{\frac{1}{3}}$ and the result follows from inequality (\ref{eqmajofftilde}).
\hfill $\square$ \vskip 2mm \noindent

\section{Other distributions and exponential concentration inequalities.}
\label{sec:other}
Using a forthcoming paper of Benaim and Rossignol
\cite{BenaimRossignolarxiv06b}, one can derive an exponential concentration
inequality on the effective resistance for various distributions, Bernoulli or
continuous ones. The only estimates needed to use the results in \cite{BenaimRossignolarxiv06b} are
the main estimates in the proof of Theorem \ref{thmmain}, and they can be
obtained easily when the resistances are bounded away from 0 and infinity. For
example, suppose that $\nu$ is bounded away
from 0 and infinity, and that it is either a Bernoulli distribution or absolutely continuous with respect to the Lebesgue measure with
a density which is bounded away from 0 on its support, then, there exist two
positive constants $C_1$ and $C_2$ such that:
$$\forall t>0,\;\PP\left(|\mathcal{R}_r(0\leftrightarrow
v)-\EE(\mathcal{R}_r(0\leftrightarrow v))|>t(\log |v|)^{\frac{1}{3}}\right)\leq C_1e^{-C_2t}\;.$$
Whether this result may be extended to distributions which are not bounded
away from 0 is still uncertain.

\section{Left-right resistance on the $n\times n$-grid.}
\label{sec:leftright}
This section is purely prospective, and focuses on another interesting case of
study: the left-right resistance on the $n\times
n$-grid on $\ZZ^2$. The graph is $\ZZ^2\cap [0,n]\times [0,n]$, the source is
$A_n=\{0\}\times [0,n]$ and the sink is
$Z_n=\{n\}\times [0,n]$. When all resistances are equal to 1, one may easily see
that $\mathcal{R}_{\overline{1}}(A_n\leftrightarrow Z_n)$ equals $n/(n+1)$ and
therefore tends to 1, as $n$ tends
to infinity. In a random setting, where all resistances are independently and
identically distributed with respect to $\nu$, this implies that:
$$\limsup_{n\rightarrow \infty}\EE(\mathcal{R}_r(A_n\leftrightarrow Z_n))\leq
\int x\;d\nu(x)\;,$$
where the inequality follows from:
$$\EE(\inf_{i\in I}f_i)\leq \inf_{i\in I} \EE(f_i)\;.$$
Recall the dual caracterisation of the resistance:
$$\frac{1}{\mathcal{R}_r(A_n\leftrightarrow Z_n)}=\inf_F\sum_{e}\frac{1}{r_e}(F(e_+)-F(e_-))^2\;,$$
where the infimum is taken over all functions $F$ on the vertices which
equal 1 on $Z_n$ and 0 on $A_n$. This implies:
$$\limsup_{n\rightarrow
  \infty}\EE\left(\frac{1}{\mathcal{R}_r(A_n\leftrightarrow Z_n)}\right)\leq
\int \frac{1}{x}\;d\nu(x)\;.$$
Using $\EE(1/R)\geq1/\EE(R)$, we finally obtain:
$$\frac{1}{\int \frac{1}{x}\;d\nu(x)}\leq\liminf_{n\rightarrow \infty}\EE(\mathcal{R}_r(A_n\leftrightarrow Z_n))\leq\limsup_{n\rightarrow \infty}\EE(\mathcal{R}_r(A_n\leftrightarrow Z_n))\leq
\int x\;d\nu(x)\;.$$
See also Theorem 2 in Hammersley \cite{Hammersley88}. In fact, it is natural
to conjecture that the limit of $\mathcal{R}_r(A_n\leftrightarrow Z_n)$ as $n$
tends to infinity exists almost surely and is constant  (see Hammersley
\cite{Hammersley88} p.~350). This is indeed the case, at least under an ellipticity condition, as follows from the work 
by K\"unnemann \cite{Kunnemann83} (see Theorem 7.4 p.~230 in the book by Jikov et al. \cite{ZhikovKozlovOleinik94}). This work relies on homogenization techniques introduced by Papanicolaou and Varadhan \cite{PapanicolaouVaradhan81}. Notice that in the book by Jikov et al., the law of large numbers is even stated for conductances which are allowed to take the value 0 (see \cite{ZhikovKozlovOleinik94}
chapters 8 and 9, notably equation
(9.16) p.~303  and Theorem 9.6, p.~314). 

Returning to resistances with finite mean, the variance in this setting
is obviously less than 1, but inequality \ref{eqmajoPoincare} suggests that it is much lower.
\begin{conj}If the resistances are bounded away from 0 and infinity,
$$Var(\mathcal{R}_r(A_n\leftrightarrow Z_n))=O\left(\frac{1}{n^2}\right)\;.$$
\end{conj}
A lower bound of this order has been proven by Wehr \cite{Wehr97} for some absolutely
continuous distributions $\nu$ under the assumption that the convergence of
the effective resistance holds almost surely. Actually, Wehr's result is
stated for effective conductivity, that is the inverse of effective
resistance, but in this context, they both are of order $\Theta(1)$, and the lower
bound of Wehr implies a lower bound of the same order on the variance of the resistance. A very appealing question is therefore:
\begin{question}
Defining
$$R_n=\mathcal{R}_r(A_n\leftrightarrow Z_n)\;,$$
does $n(R_n-\EE(R_n))$ converge in distribution as $n$ tends to infinity ?
What is the limit law ?
\end{question}



\section{Submean variance bound for $p$-resistance}
\label{sec:penergy}
In the analysis presented so far, the probabilistic interpretation of
electrical networks has played no role. It is therefore tempting to extend our
work to the setting of \emph{$p$-networks} (see for instance Soardi
\cite{Soardi94} p.176-178). As before, consider an unoriented, at most
countable and locally finite graph $G=(V,E)$. Let
$r=(r_e)_{e\in E}$ be a collection of resistances. For any $p>1$, we define the
$p$-resistance between two vertices $x$ and $y$ as
\begin{equation}
\label{eqdefpresistance}
\mathcal{R}^p_r(x\leftrightarrow y)=\inf_{\|\theta\|=1}\sum_{e\in E}
r_e|\theta(e)|^p\;,
\end{equation}
where the infimum is taken over all flows $\theta$ from $x$ to $y$ with
strength 1. It is known that the $p$-resistance from 0 to infinity on $\ZZ^d$, when all
resistances equal 1, is finite
if and only if $p> d/(d-1)$ (see Soardi and Yamasaki
\cite{SoardiYamasaki93}). More precisely, the flow described in
\cite{SoardiYamasaki93} and the usual shorting argument to lowerbound
resistance from 0 to the border of the box $B_n$ allow to obtain the following
estimate of the $p$-resistance on $\ZZ^d$:
$$\EE(\mathcal{R}^p_r(0\leftrightarrow
v))=\Theta\left(\sum_{k=0}^{|v|}\frac{1}{(2n+1)^{(d-1)(p-1)}}\right)\;.$$
Whenever this expectation tends to infinity as $|v|$ tends to infinity, and when
$d\geq 2$, one may hope
to obtain a similar result as in Theorem \ref{thmmain}. This is indeed the
case: we obtain a weaker result, but still, the variance is small compared to
the mean. The proof is essentially the same as in the case where $p=2$. There are
two main important points to take care of. First, it remains true that for a unit flow which
minimizes the $p$-energy $\sum_{e\in E}
r_e|\theta(e)|^p$, the flow on each edge is less than 1. This follows from the same argument as in the linear case (see Lyons and
Peres \cite{LyonsPeresinprogress}, p.~49-50). Second, it is not clear whether
inequality (\ref{eqpseudotriangular}) remains true or not. Nevertheless, we can easily obtain the following weaker inequality. For every
three vertices $x,y,z$ in $\ZZ^d$, and any $r\in\Omega$,
\begin{equation}
\label{eqpseudotriangularpnetwork}
\mathcal{R}^p_r(x,z)\leq \mathcal{R}^p_r(x,y)+2^pb|z-y|\;.
\end{equation}
To see this, let $\theta^{x,y}$ be
the unit current flow (for $p$-energy) from $x$ to $y$ and $\pi=(u_0=y,u_1,\ldots , u_{|z-y|}=z)$ be a deterministic oriented path from $y$
to $z$. Define a flow $\theta^{y,z}$ from $y$ to $z$ as
follows:
$$
\begin{array}{rcl}\theta^{y,z}(e)&=&0,\;\mbox{ if }e\not\in \pi\;,\\
\theta^{y,z}(\overrightarrow{u_iu_{i+1}})&=&1,\;\mbox{ if
}i\in\{0,\ldots,|z-y|-1\}\end{array}\;.
$$
Now, let $\eta^{x,z}$ be the unit flow $\theta^{x,y}+\theta^{y,z}$, which goes from
$x$ to $z$.
\begin{eqnarray*}
\mathcal{R}^p_r(x,z)&\leq &\sum_{e}r_e|\eta^{x,z}(e)|^p\;,\\
&=&\sum_{e\not\in\pi}r_e|\theta^{x,y}(e)|^p+\sum_{e\in\pi}r_e|\theta^{x,y}(e)+1|^p\;,\\
&\leq&\mathcal{R}^p_r(x,y)+2^pb|z-y|\;,
\end{eqnarray*}
where the last inequality follows from the fact that $|\theta^{x,y}_e|\leq 1$
for every edge $e$. Inequality (\ref{eqpseudotriangularpnetwork}) is
proved. This allows to adapt the proof of Theorem \ref{thmmain} as follows. Inequalities (\ref{eqmajofftilde}) and
(\ref{eqmajovarz}) become respectively:
$$Var(f_v)\leq (Cm)^2+2Cm\sqrt{Var(\tilde f)}+Var(\tilde f)\;,$$
and
$$\EE_\nu(Var_\mu(\tilde f))=O(m^2)\;.$$
The rest of the proof is the same, and leads to:
$$Var (\tilde f)= O\left(\frac{\EE(f)}{\log m}+m^2\right)\;.$$
We can choose, for instance $m=\lceil (\EE(f))^{1/3}\rceil$, to get the
following weaker analog of Theorem \ref{thmmain}.

\begin{prop}
\label{proppenergy}
Suppose that $\nu=\frac{1}{2}\delta_a+\frac{1}{2}\delta_b$, with $0<a\leq
b<+\infty$. Let $d\geq 1$ be an integer, $E$ be the set of edges in $\ZZ^d$, and define
$\mu=\nu^{E}$. Then, for any real number $p$ in $]1,+\infty[$:
$$\EE(\mathcal{R}^p_r(0\leftrightarrow
v))=\Theta(a_d(|v|,p))\;,$$
and if $d\geq 2$,
$$Var_\mu(\mathcal{R}^p_r(0\leftrightarrow
v))=O\left(\frac{a_d(|v|,p)}{\log a_d(|v|,p)}\right)\;,$$
where
$$a_d(n,p)=\left\lbrace\begin{array}{ll}n^{1-(d-1)(p-1)}&\mbox{ if }p<\frac{d}{d-1}\;,\\
\log (n)&\mbox{ if }p=\frac{d}{d-1},\\
1&\mbox{ if }p>\frac{d}{d-1}.\end{array}\right.
$$
\end{prop} 
Remark also that Lemma \ref{lemminovariance} is easily extended to this setting, and we get
therefore that, for any $p>\frac{d}{d-1}$,
$$Var_\mu(\mathcal{R}^p_r(0\leftrightarrow v))=\Theta(1)\;.$$

\section*{Acknowledgements}
Itai Benjamini would like to thank Gady Kozma, Noam Berger and Oded Schramm for
useful discussions. R. Rossignol would like to thank Michel Benaim for useful
discussions, and Yuval Peres for useful remarks on a first version of this
paper.



\begin{thebibliography}{}

\bibitem[An{\'e} et~al., 2000]{thesardsledoux}
An{\'e}, C., Blach{\`e}re, S., Chafa{\"{\i}}, D., Foug{\`e}res, P., Gentil, I.,
  Malrieu, F., Roberto, C., and Scheffer, G. (2000).
\newblock {\em Sur les in{\'e}galit{\'e}s de Sobolev logarithmiques}.
\newblock Soci{\'e}t{\'e} Math{\'e}matique de France, Paris.

\bibitem[Benaim and Rossignol, 2006a]{BenaimRossignolarxiv06b}
Benaim, M. and Rossignol, R. (2006a).
\newblock Exponential concentration for {F}irst {P}assage {P}ercolation through
  modified {P}oincar\'e inequalities.
\newblock \verb+http://arxiv.org/abs/math.PR/0609730+ to appear in {A}nnales de
  l'{IHP}.

\bibitem[Benaim and Rossignol, 2006b]{BenaimRossignolarxiv06}
Benaim, M. and Rossignol, R. (2006b).
\newblock A modified {P}oincar\'e inequality and its application to first
  passage percolation.
\newblock \verb+http://arxiv.org/abs/math.PR/0602496+.

\bibitem[Benjamini et~al., 2003]{BenjaminiKalaiSchramm03}
Benjamini, I., Kalai, G., and Schramm, O. (2003).
\newblock First passage percolation has sublinear distance variance.
\newblock {\em Ann. Probab.}, 31(4):1970--1978.

\bibitem[Berger, 2002]{Berger02}
Berger, N. (2002).
\newblock Transience, recurrence and critical behavior for long-range
  percolation.
\newblock {\em Commun. Math. Phys.}, 226:531--558.

\bibitem[Doyle and Snell, 1984]{DoyleSnell84}
Doyle, P. and Snell, J. (1984).
\newblock {\em Random walks and electric networks}.
\newblock Mathematical {A}ssociation of {A}merica.
\newblock Also available at the arxiv as math.PR/0001057.

\bibitem[Falik and Samorodnitsky, 2006]{FalikSamorodnitsky06}
Falik, D. and Samorodnitsky, A. (2006).
\newblock Edge-isoperimetric inequalities and influences.
\newblock to appear \verb+http://arxiv.org/pdf/math.CO/0512636+.

\bibitem[Hammersley, 1988]{Hammersley88}
Hammersley, J.~M. (1988).
\newblock Mesoadditive processes and the specific conductivity of lattices.
\newblock {\em J. Appl. Probab.}, Special Vol. 25A:347--358.

\bibitem[Jikov et~al., 1994]{ZhikovKozlovOleinik94}
Jikov, V.~V., Kozlov, S.~M., and Oleinik, O.~A. (1994).
\newblock {\em Homogenization of differential operators and integral
  functionals}.
\newblock Springer-Verlag.

\bibitem[Kesten, 1986]{Kesten84}
Kesten, H. (1986).
\newblock Aspects of first passage percolation.
\newblock In {\em Ecole d'{\'et\'e} de probabilit{\'e} de Saint-Flour
  XIV---1984}, volume 1180 of {\em Lecture Notes in Math.}, pages 125--264.
  Springer, Berlin.

\bibitem[K{\"u}nnemann, 1983]{Kunnemann83}
K{\"u}nnemann, R. (1983).
\newblock The diffusion limit for reversible jump processes on
  $\mathbb{Z}^{d}$\ with ergodic random bond conductivities.
\newblock {\em Comm. Math. Phys.}, 90(1):27--68.

\bibitem[Lyons and Peres, 2006]{LyonsPeresinprogress}
Lyons, R. and Peres, Y. (1997-2006).
\newblock {\em Probability on trees and networks}.
\newblock Book in progress.
\newblock Available at
  \verb+http://mypage.iu.edu/~rdlyons/prbtree/prbtree.html+.

\bibitem[Papanicolaou and Varadhan, 1981]{PapanicolaouVaradhan81}
Papanicolaou, G. and Varadhan, S. (1981).
\newblock Boundary value problems with rapidly oscillating random coefficients.
\newblock In {\em Random fields, Vol. I, II (Esztergom, 1979)}, volume~27 of
  {\em Colloq. Math. Soc. J\'anos Bolyai}, pages 835--873. North-Holland,
  Amsterdam.

\bibitem[Pemantle and Peres, 1996]{PemantlePeres96}
Pemantle, R. and Peres, Y. (1996).
\newblock On which graphs are all random walks in random environments transient
  ?
\newblock {\em Random Disc. Struct.}, 76:207--211.

\bibitem[Peres, 1999]{Peres97}
Peres, Y. (1999).
\newblock Probability on trees: an introductory climb.
\newblock In {\em Lectures on probability theory and statistics (Saint-Flour,
  1997)}, volume 1717, pages 193--280, Berlin. Springer.

\bibitem[Rossignol, 2006]{Rossignol06}
Rossignol, R. (2006).
\newblock Threshold for monotone symmetric properties through a logarithmic
  {S}obolev inequality.
\newblock {\em Ann. Probab.}, 34(5):1707--1725.

\bibitem[Soardi, 1994]{Soardi94}
Soardi, P. (1994).
\newblock {\em Potential theory on infinite networks}.
\newblock Number 1590 in Lecture Notes in Mathematics. Springer-Verlag, Berlin.

\bibitem[Soardi and Yamasaki, 1993]{SoardiYamasaki93}
Soardi, P. and Yamasaki, M. (1993).
\newblock Parabolic index and rough isometries.
\newblock {\em Hiroshima Math. J.}, 23:333--342.

\bibitem[Steele, 1986]{Steele86}
Steele, J. (1986).
\newblock An {E}fron-{S}tein inequality for nonsymmetric statistics.
\newblock {\em Ann. Stat.}, 14:753--758.

\bibitem[Talagrand, 1994]{Talagrand94a}
Talagrand, M. (1994).
\newblock On {R}usso's approximate zero-one law.
\newblock {\em Ann. Probab.}, 22:1576--1587.

\bibitem[Wehr, 1997]{Wehr97}
Wehr, J. (1997).
\newblock A lower bound on the variance of conductance in random resistor
  networks.
\newblock {\em J. Statist. Phys.}, 86(5-6):1359--1365.

\end{thebibliography}


\end{document}